\documentclass[twocolum,fleqn]{autart}    

\usepackage{graphicx}          

\usepackage{amssymb}
\usepackage[cmex10]{amsmath}
\usepackage{accents}
\usepackage{epsfig}

\usepackage{bm}
\usepackage{amsmath}
\usepackage{mathptmx}
\usepackage{color}
\usepackage{dsfont}
\usepackage{empheq}
\usepackage{hyperref}
\usepackage{array}
\usepackage{cite}
\allowdisplaybreaks
\usepackage{subfig}
\usepackage{pgfplots}
\usepackage{tikz}
\usetikzlibrary{arrows,shapes,patterns,positioning,calc}
\tikzstyle{boxstyle}=[draw=black,inner sep=7pt]     
\tikzstyle{arrowstyle}=[]    
\usetikzlibrary{arrows, decorations.markings}

\newlength\figureheight
\newlength\figurewidth

\newtheorem{remark}{Remark}
\newtheorem{lemma}{Lemma}
\newtheorem{theorem}{Theorem}
\newcommand{\pa}{\partial}
\newcommand{\ba}{\begin{align}}
\newcommand{\ea}{\end{align}}
\newcommand{\fr}{\frac}
\newcommand{\gm}{\gamma}
\newcommand{\ep}{\varepsilon}

\begin{document}

\begin{frontmatter}

\title{Sampled-Data Control of the Stefan System} 
 

\author[UCSD]{Shumon Koga}\ead{skoga@eng.ucsd.edu},               
\author[iasson]{Iasson Karafyllis}\ead{iasonkar@central.ntua.gr}, 
\author[UCSD]{Miroslav Krstic}\ead{krstic@ucsd.edu}  

\address[UCSD]{Department of Mechanical and Aerospace Engineering, University of California, San Diego, La Jolla, CA 92093-0411 USA}                                               
\address[iasson]{Department of Mathematics, National Technical University of Athens, 15780 Athens, Greece}
          

\begin{abstract}                          
This paper presents results for the sampled-data boundary feedback control to the Stefan problem. The Stefan problem represents a liquid-solid phase change phenomenon which describes the time evolution of a material's temperature profile and the interface position. First, we consider the sampled-data control for the one-phase Stefan problem by assuming that the solid phase temperature is maintained at the equilibrium melting temperature. We apply Zero-Order-Hold (ZOH) to the nominal continuous-time control law developed in \cite{Shumon19journal} which is designed to drive the liquid-solid interface position to a desired setpoint. Provided that the control gain is bounded by the inverse of the upper diameter of the sampling schedule, we prove that the closed-loop system under the sampled-data control law satisfies some conditions required to validate the physical model, and the system's origin is globally exponentially stable in the spatial $L_2$ norm. Analogous results for the two-phase Stefan problem which incorporates the dynamics of both liquid and solid phases with moving interface position are obtained by applying the proposed procedure to the nominal control law for the two-phase problem developed in \cite{koga2018CDC}. Numerical simulation illustrates the desired performance of the control law implemented to vary at each sampling time and keep constant during the period. 
\end{abstract}

\end{frontmatter}

\section{Introduction} 
\subsection{Background} 
Liquid-solid phase transitions are physical phenomena which appear in various kinds of science and engineering processes. Representative applications include sea-ice melting and freezing~\cite{koga2019arctic}, continuous casting of steel \cite{petrus2012}, cancer treatment by cryosurgeries \cite{Rabin1998}, additive manufacturing for materials of both polymer \cite{koga2018polymer} and metal \cite{chung2004}, crystal growth~\cite{conrad_90}, lithium-ion batteries \cite{koga2017battery}, and thermal energy storage systems \cite{zalba03}. Physically, these processes are described by a temperature profile along a liquid-solid material, where the dynamics of the liquid-solid interface is influenced by the heat flux induced by melting or solidification. A mathematical model of such a physical process is called the Stefan problem\cite{Gupta03}, which is formulated by a diffusion PDE defined on a time-varying spatial domain. The domain's length dynamics is described by an ODE dependent on the gradient of the PDE state. Apart from the thermodynamical model, the Stefan problem has been employed to model several chemical, electrical, social, and financial dynamics such as tumor growth process \cite{Friedman1999}, domain walls in ferroelectric thin films \cite{mcgilly2015}, spreading of invasive species in ecology \cite{Du2010speading}, information diffusion on social networks \cite{Lei2013}, and optimal exercise boundary of the American put option on a zero dividend asset \cite{Chen2008}.

 While the numerical analysis of  the one-phase Stefan problem is broadly covered in the literature, their control related problems have been rarely addressed. In addition to it, most of the proposed control approaches are based on finite dimensional approximations with the assumption of  an explicitly given moving boundary dynamics \cite{Daraoui2010,Armaou01}. For control objectives, infinite-dimensional approaches have been used for stabilization of  the temperature profile and the moving interface of a 1D Stefan problem, such as enthalpy-based feedback~\cite{petrus2012} and geometric control~\cite{maidi2014}.  These works designed control laws ensuring the asymptotical stability of the closed-loop system in the ${L}_2$ norm. However, the results in \cite{maidi2014} are established based on the assumptions on the liquid temperature being greater than the melting temperature, which must be ensured by showing the positivity of the boundary heat input. 
 

Recently, boundary feedback controllers for the Stefan problem have been designed via a ``backstepping transformation" \cite{krstic2008boundary,andrew2004} which has been used for many other classes of infinite-dimensional systems. For instance, \cite{Shumon16} designed a state feedback control law by introducing a nonlinear backstepping transformation for moving boundary PDE, which achieved the exponentially stabilization of the closed-loop system in the ${\mathcal H}_1$ norm without imposing any {\em a priori} assumption. Based on the technique, \cite{Shumon16CDC} designed an observer-based output feedback control law for the Stefan problem, \cite{Shumon19journal} extended the results in \cite{Shumon16, Shumon16CDC} by studying the robustness with respect to the physical parameters and developed an analogous design with Dirichlet boundary actuation, \cite{Shumon17ACC} designed a state feedback control for the Stefan problem under the material's convection, \cite{koga_2019delay} developed a control design with time-delay in the actuator and proved a delay-robustness, \cite{koga2019iss} investigated an input-to-state stability of the control of Stefan problem with respect to an unknown heat loss at the interface, and \cite{koga2018CDC} developed a control design for the two-phase Stefan problem. 

The aforementioned results assumed the control input to be varying continuously in time; however, in practical implementation of the control systems it is impossible to dynamically change the control input continuously in time due to limitations of the sensors, actuators, and software. Instead, the control input can be adjusted at each sampling time at which the measured states are obtained or the actuator is manipulated. One of the most fundamental and well known method to design such a ``sampled-data" control is the so-called ``emulation design" that applies ``Zero-Order-Hold" (ZOH) to the nominal ``continuous-time" control law. A general result for nonlinear ODEs to guarantee the global stability of the closed-loop system under such a ZOH-based sampled-data control was studied in \cite{Karafyllis09}, and the sampled-data observer design under discrete-time measurement is developed in \cite{Karafyllis09observer} by introducing inter-sampled output predictor. As further extensions, the stability of the sampled-data control for general nonlinear ODEs under actuator delay is shown in \cite{Karafyllis12,karafyllis2017book} by applying predictor-based feedback developed in \cite{krstic2009delay}, and results for a linear parabolic PDE are given in \cite{Karafyllis18sampledpde} by employing Sturm-Liouville operator theory. The sampled-data control for parabolic PDEs has been intensively developed by Fridman and coworkers by utilizing linear matrix inequalities \cite{Am14, Emilia12,Emilia13, Selivanov16}. However, none of the existing work on the sampled-data control has studied the class of the Stefan problem described by a parabolic PDE with state-dependent moving boundaries ``(a nonlinear system)". 

\subsection{Contributions and results} 
This paper presents the first theoretical result for the sampled-data boundary feedback control for the Stefan problem. The approach employed in this paper is distinct from the methodology developed in literature. Namely, we solve the growth of the system's energy analytically in time under the proposed sampled-data feedback control that is in the form of an energy-shaping design. Then, a perturbation that is incorporated in the closed-loop system due to the error between the continuous-time design and the sampled-data design can be represented analytically, and the closed-loop stability is proven by using Lyapunov method. 

First, we consider the one-phase Stefan problem by assuming that the solid phase temperature is maintained at the melting temperature and focusing on the single melting process. We employ ZOH to the nominal continuous-time feedback controller for the one-phase Stefan problem developed in \cite{Shumon19journal}, and prove the required conditions for the model validity and the global exponential stability of the closed-loop system under explicit conditions for the setpoint position and the control gain with respect to the sampling scheduling. Next, we consider the two-phase Stefan problem by incorporating the dynamics of the solid phase temperature and prove the analogous results for the sampled-data control for the two-phase Stefan problem. The results established in this paper hold for arbitrary sampling schedules, and not necessarily uniform sampling schedules. 

\subsection{Organization} 
The mathematical model the one-phase Stefan problem for a single phase change is presented in Section \ref{sec:model} with stating some important properties. The sampled-data control law and the stability proof of the closed-loop system is given in Section \ref{sec:control}. The extension of the presented procedure to the two-phase Stefan problem is described in Section \ref{sec:twophase}. The numerical simulation of the proposed control law is provided in Section \ref{sec:simulation}. The paper ends with the concluding remarks in Section \ref{sec:conclusion}.

\section{Description of the One-Phase Stefan Problem} \label{sec:model} 
\begin{figure}[t]
\centering
\includegraphics[width=2.5in]{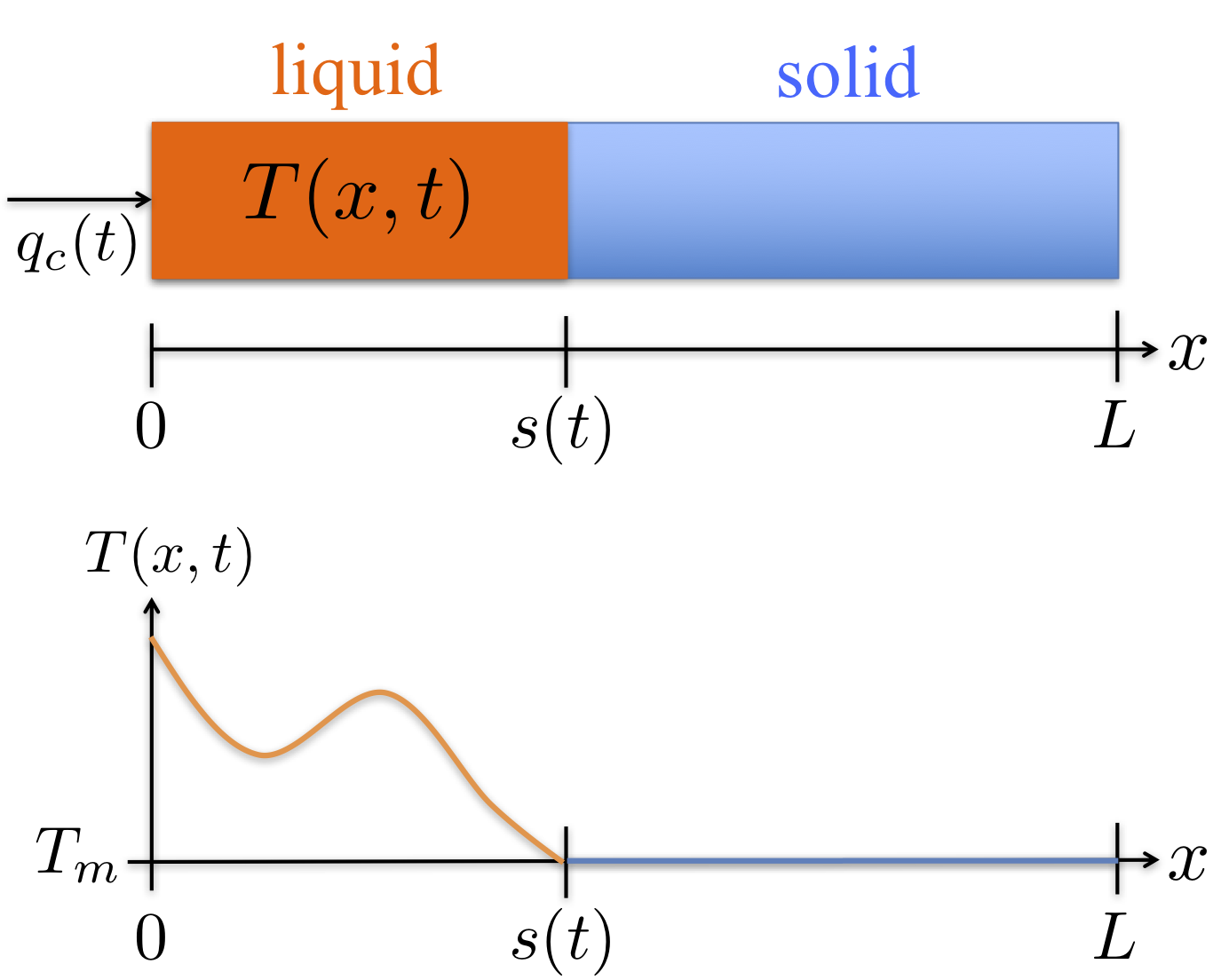}\\
\caption{Schematic of one-phase Stefan problem.}
\label{fig:stefan}
\end{figure}

Consider a physical model which describes the melting or solidification mechanism in a pure one-component material of length $L$ in one dimension. In order to mathematically describe the position at which phase transition occurs, we divide the domain $[0, L]$ into two time-varying sub-domains, namely,  the interval $[0,s(t)]$ which contains the liquid phase, and the interval $[s(t),L]$ that contains the solid  phase. A heat flux enters the material through the boundary at $x=0$ (the fixed boundary of the liquid phase) which affects the liquid-solid interface dynamics through heat propagation in liquid phase. As a consequence, the heat equation alone does not provide a complete description of the phase transition and must be coupled with the dynamics that describes the moving boundary. This configuration is shown in Fig.~\ref{fig:stefan}.

Assuming that the temperature in the liquid phase is not lower than the melting temperature  of the material $T_{{\mathrm m}}$, the energy conservation and heat conduction laws yield the heat equation of  the liquid phase as follows
\begin{align}\label{sys1}
T_t(x,t)&=\alpha T_{xx}(x,t), \quad  \alpha :=\fr{k}{\rho C_{{\rm p}}}, \quad  0\leq x\leq s(t), 
\end{align}
with the boundary conditions
\begin{align}\label{sys2}
-k T_x(0,t)&=q_{\rm c}(t), \\ \label{sys3}
T(s(t),t)&=T_{{\rm m}},
\end{align}
and the initial values
\begin{align}
T(x,0)=T_0(x), \quad s(0) = s_0
\end{align}
where $T(x,t)$,  ${q}_{\rm c}(t)$,  $\rho$, $C_{{\rm p}}$, and $k$ are the distributed temperature of the liquid phase, the manipulated heat flux, the liquid density, the liquid heat capacity, and the liquid heat conductivity, respectively. Moreover, the local energy balance at the liquid-solid interface $x=s(t)$ yields
\begin{align}\label{sys4}
 \dot{s}(t)=- \beta T_x(s(t),t), \quad \beta := \frac{k}{\rho \Delta H^*} 
\end{align}
where $\Delta H^*$ represents the latent heat of fusion. 
\begin{remark}As the moving interface  $s(t)$ depends on the temperature, the problem defined in  \eqref{sys1}--\eqref{sys4}  is nonlinear.\end{remark}

There are two underlying assumptions to validate the model \eqref{sys1}-\eqref{sys4}. First, the liquid phase is not frozen to the solid phase from the boundary $x=0$. This condition is ensured if the liquid temperature $T(x,t)$ is greater than the melting temperature. Second, the material is not completely melt or frozen to single phase through the disappearance of the other phase. This condition is guaranteed if the interface position remains inside the material's domain. In addition, these conditions are also required for the well-posedness (existence and uniqueness) of the solution in this model. Taking into account of these model validity conditions, we emphasize the following remark. 
\begin{remark} \emph{
To maintain the model \eqref{sys1}-\eqref{sys4} to be physically validated, the following conditions must hold: 
\begin{align}\label{temp-valid}
T(x,t) \geq& T_{{\rm m}}, \quad  \forall x\in(0,s(t)), \quad \forall t>0, \\
\label{int-valid}0<&s(t)<L, \quad \forall t>0. 
\end{align}}
\end{remark}

Based on the above conditions, we impose the following assumption on the initial data. 
\begin{assum}\label{initial} 
$s_0 \in (0,L)$, $T_0(x) \geq T_{{\rm m}}$ for all $x \in [0,s_0]$, and $T_0(x) $ is continuously differentiable in $x\in[0,s_0]$.
 \end{assum}
 The existence and uniqueness of the solution of the one-phase Stefan problem \eqref{sys1}--\eqref{sys4} is presented in \cite{Cannon1971} as follows. 
\begin{lemma}\label{lem1}
Under Assumption \ref{initial}, if $q_{\rm c}(t)$ is a bounded piecewise continuous function with generating nonnegative heat for a time interval, i.e., $q_{{\rm c}}(t) \geq 0, $ for all $t \in [0, \bar t]$, then there exists a unique solution for the Stefan problem \eqref{sys1}-\eqref{sys4} with satisfying the condition \eqref{temp-valid} for all $t \in [0, \bar t]$. Moreover, it holds 
\begin{align} 
\dot{s}(t)>0, \quad \forall t \in [0, \bar t]. 	
\end{align}
\end{lemma}

\section{Sampled-Data Control for the One-Phase Stefan Problem} \label{sec:control} 

\subsection{Problem statement and main result} 
The steady-state solution $(T_{{\rm eq}}(x), s_{{\rm eq}})$ of the system \eqref{sys1}--\eqref{sys4} with zero manipulating heat flux $q_{\rm c}(t)=0$ yields a uniform melting temperature $T_{{\rm eq}}(x) = T_{{\mathrm m}}$ and a constant interface position given by the initial data. In \cite{Shumon16}, the authors developed the exponential stabilization of the interface position $s(t)$ at a desired reference setpoint $s_{{\mathrm r}}$ through the design of $q_{\rm c}(t)$ as
\begin{align} \label{nominalcont} 
q_{\rm c}(t) = - c \left(\frac{k}{\alpha} 	\int_0^{s(t)} (T(x,t) - T_{{\rm m}}) {\rm d}x + \frac{k}{\beta} (s(t) - s_r ) \right) ,  
\end{align}
where $c>0$ is the controller gain. However, in practical implementation, the actuation value cannot be changed continuously in time. Instead, by obtaining the measured value as signals discretely in time, the control value needs to be implemented at each sampling time. One of the most typical design for such a sampled-data control is the application of "Zero-Order-Hold"(ZOH) to the nominal continuous time control law. Through ZOH, during the time intervals between each sampling, the control maintains the value at the previous sampling time. Let $t_j$ be the $j$-th sampling time for $j = 0, 1, 2, \cdots, $, and $\tau_j$ be defined by 
\begin{align}
\tau_{j} = t_{j+1} - t_j . 
\end{align}  
The application of ZOH to the nominal control law \eqref{nominalcont} leads to the following design for the sampled-data control 
\begin{align} \label{sampledcont}
q_{\rm c}(t) = &- c \left(\frac{k}{\alpha} 	\int_0^{s(t_j)} (T(x,t_j) - T_{{\rm m}}) {\rm d}x + \frac{k}{\beta} (s(t_j) - s_r ) \right), \notag\\
& \forall t \in [t_{j}, t_{j+1}),  
\end{align}
 of which the right hand side is constant during the time interval $t\in [t_{j}, t_{j+1}) $. Let us denote $q_{j} = q_{\rm c}(t)$ for $t \in[t_{j}, t_{j+1} ) $. Hereafter, all the variables with subscript $j$ denote the variables at $t = t_{j}$. First, we introduce the following assumptions on the setpoint and the sampling scheduling. 
 \begin{assum}  \label{ass1}
The setpoint is chosen to verify
\begin{align}\label{eq:setpoint}
 s_0 + \fr{\beta}{\alpha} \int_0^{s_0} (T_0(x)-T_{\rm m}) {\rm d}x < s_r < L. 
\end{align}
\end{assum}

\begin{assum} \label{ass:sample} 
The sampling schedule has a finite upper diameter and a positive lower diameter, i.e., there exist constants $0 < r \leq R$ such that 
\begin{align} 
\sup_{j \in {\mathcal Z}^+}\{\tau_j\} \leq & R , \\
\inf_{j \in {\mathcal Z}^+} \{\tau_j\} \geq & r. 
\end{align} 
\end{assum} 



Our main theorem is given next. 
\begin{theorem}\label{theo-1}
Consider the closed-loop system \eqref{sys1}--\eqref{sys3}, \eqref{sys4}, \eqref{sampledcont} under Assumptions \ref{initial}, \ref{ass1}. Then for every $0<r \leq R<1/c$, there exists a constant $M:=M(r)$ for which the following property holds: for every sequence $\{t_j\geq 0 :  j = 0, 1, 2, \dots \}$ with $t_0=0$ for which Assumption 3 holds, the initial-boundary value problem \eqref{sys1}--\eqref{sys4} with \eqref{sampledcont} has a unique solution satisfying \eqref{temp-valid}, \eqref{int-valid} as well as the following estimate:
\begin{align}\label{h1}
\Psi(t) \leq M \Psi(0) \exp(-bt) ,
\end{align}
where $b=\fr{1}{8} \min \left\{ \fr{\alpha}{s_r^2}, c \right\}$, for all $t\geq0$, in the $L_2$ norm $\Psi(t) =  \int_{0}^{s(t)} \left(T(x,t)-T_{\rm m}\right)^2 {\rm d}x +(s(t)-s_{r})^2$.
\end{theorem}

The proof of Theorem \ref{theo-1} is established through several steps in the next sections. The positive constant $M$ in \eqref{h1} has a dependency on $r>0$ as 
\begin{align} 
M(r) = M_1 + \fr{M_2}{1 - \left( 1 - c r \right)^2   e^{\frac{ cr}{8}} } , 
\end{align} 
for some positive constants $M_1>0$ and $M_2>0$ that are not dependent on $r>0$. 

\subsection{Some key properties of the closed-loop system}
We first provide the following lemma. 
\begin{lemma} \label{lem2} 
The closed-loop system consisting of the plant \eqref{sys1}--\eqref{sys4} under the sampled-data control law \eqref{sampledcont}	has a unique classical solution which is equivalent to the open-loop solution of \eqref{sys1}--\eqref{sys4} with the control law of 
\begin{align} \label{qjsol}
	q_{\rm c}(t)  = q_{j}=  & q_{0}\prod_{i=0}^{j-1} \left( 1 - c \tau_i \right), \quad \forall t \in [t_{j}, t_{j+1}), \quad \forall j \in {\mathcal Z}^{+}
\end{align}
where 
\begin{align} 
q_0 = 	- c \left(\frac{k}{\alpha} 	\int_0^{s_0} (T_0(x) - T_{{\rm m}}) {\rm d}x + \frac{k}{\beta} (s_0 - s_r ) \right).  
\end{align}

\end{lemma}
\begin{pf} 
We introduce the following reference error states:  
\begin{align}
u(x,t) = T(x,t) - T_{\rm m}, \quad X(t) = s(t) - s_r. 
\end{align}
The governing equations \eqref{sys1}--\eqref{sys4} are rewritten as the following reference error system 
\begin{align}\label{u-sys1}
u_{t}(x,t) =&\alpha u_{xx}(x,t),\\
\label{u-sys2}u_x(0,t) =& - q_{\rm c}(t)/k,\\
\label{u-sys3}u(s(t),t) =&0,\\
\label{u-sys4}\dot{X}(t) =&-\beta u_x(s(t),t) .
\end{align}
Define the internal energy of the reference error system as follows:  
\begin{align} \label{Energy}
\widetilde E(t) =  \fr{k}{\alpha } \int_0^{s(t)} u(x,t) {\rm d}x + \fr{k}{\beta} X(t). 
\end{align}
Taking the time derivative of \eqref{Energy} along the solution of \eqref{u-sys1}--\eqref{u-sys4} leads to the following energy conservation law
\begin{align} \label{conservation}
\fr{{\rm d}}{{\rm d}t} \widetilde E(t) = q_{\rm c}(t). 
\end{align}
Noting that $q_{\rm c}(t) $ is constant for $t \in [t_{j}, t_{j+1})$ as $q_{\rm c}(t) = q_{j}$ under ZOH-based sampled-data control, taking the integration of \eqref{conservation} from $t = t_{j}$ to $t = t_{j+1}$ yields 
\begin{align} \label{Ejqj}
\widetilde E_{j+1} - \widetilde E_{j} = \tau_j q_{j} , 
\end{align}
where $\widetilde E_{j} = \widetilde E(t_{j})$ and $\tau_{j} = t_{j+1} - t_{j}$. The sampled-data control \eqref{sampledcont} and the internal energy \eqref{Energy} at each sampling time satisfy the following relation: 
\begin{align} \label{qj}
q_{j} = - c \widetilde E_{j}. 
\end{align}
Substituting \eqref{qj} into \eqref{Ejqj}, we obtain
\begin{align}
\widetilde E_{j+1} = \left(1 - c \tau_j \right) \widetilde E_{j}, 
\end{align}
which leads to the explicit solution as follows:  
\begin{align}\label{Ejsol} 
\widetilde E_{j} =& \widetilde E_0 \prod_{i=0}^{j-1} \left( 1 - c \tau_i \right) . 
\end{align}
Substituting \eqref{Ejsol} into \eqref{qj} yields as \eqref{qjsol}. 
Therefore, the closed-loop system under the sampled-data feedback control \eqref{sampledcont} is equivalent to the open-loop solution with the control input \eqref{qjsol}. Moreover, under Assumptions \ref{ass1}, \ref{ass:sample}, and the fact that $c<\fr{1}{R}$, the input \eqref{qjsol} is shown to be a bounded piecewise continuous function and $q_{\rm c}(t) \geq 0$ for all $t \geq 0$. Thus, the existence and uniqueness of the solution is ensured by Lemma \ref{lem1}, from which we conclude Lemma \ref{lem2}. 
\end{pf} 

\begin{lemma} \label{lem3} 
The closed-loop system satisfies the following properties: 
\begin{align} \label{prop1}
\dot s(t) > 0, \quad \forall t \geq 0, \\
s_0 < s(t) < s_r, \quad \forall t \geq 0. 	\label{prop2}
\end{align}
\end{lemma}
\begin{pf} 
Combining Lemma \ref{lem1} with Lemma \ref{lem2}, one can deduce \eqref{prop1}, and $s_0<s(t)$ for all $t \geq 0$. We show $s(t)<s_r$ for all $t \geq 0$. Integrating \eqref{conservation} from $t = t_{j}$ to $t \in [t_{j}, t_{j+1})$ leads to 
\begin{align}\label{Etqj} 
\widetilde E(t) - \widetilde E_{j} = (t - t_{j} ) q_{j}, \quad \forall t \in [t_j, t_{j+1}). 
\end{align}
With the help of \eqref{qj} and \eqref{Ejsol}, equation \eqref{Etqj} yields
\begin{align} \label{EtEk}
\widetilde E(t) = \left( 1 - c (t - t_{j} )  \right) \widetilde  E_{j}, \quad \forall t \in [t_{j}, t_{j+1}).  
\end{align}
By Assumption \ref{ass:sample} and since $c<\fr{1}{R}$, we have $0<c < \frac{1}{\tau_{j}}$ for all $j \in {\mathcal Z}^+$. In addition, for all $ t \in [t_{j}, t_{j+1})$ and for all $j \in {\mathcal Z}^+$, it holds $t - t_{j} \leq \tau_{j} $. Hence, we have $ 1 - c (t - t_{j} ) >0 $, for all $ t \in [t_{j}, t_{j+1} )$ and for all $j \in {\mathcal Z}^+$. Applying this to \eqref{EtEk} and noting that 
\begin{align} 
\widetilde E_j < 0, \quad \forall j \in {\mathcal Z}^+, 	
\end{align}
deduced from \eqref{Ejsol} and Assumption \ref{ass1}, one can obtain 
\begin{align} \label{Etnegative} 
\widetilde E(t) < 0, \quad \forall t \geq0. 	
\end{align}
Substituting \eqref{Etnegative} into \eqref{Energy} and applying $u(x,t) >0$ for all $x \in (0,s(t))$ and $t \geq 0$, we have
\begin{align} \label{Xtnegative} 
X(t) < 0, \quad \forall t \geq0	, 
\end{align}
which leads to $s(t) < s_r$ for all $t \geq 0$.
\end{pf} 

  \subsection{Stability analysis} \label{sec:stability}
  To conclude Theorem \ref{theo-1}, this section is devoted to the stability proof of the closed-loop system under the designed sampled-data control law. First, we introduce the backstepping transformation developed in \cite{Shumon19journal} for the continuous-time design, and apply the transformation to the closed-loop system under the sampled-data control in this paper. 
  
\subsubsection{State transformation}
  
  Introduce the following backstepping transformation
\begin{align}\label{eq:DBST}
w(x,t)=&u(x,t)-\frac{\beta}{\alpha} \int_{x}^{s(t)} \phi (x-y)u(y,t) {\rm d}y \notag\\
&-\phi(x-s(t)) X(t),
\end{align}
which maps into
\begin{align}\label{tarPDE}
w_t(x,t)=&\alpha w_{xx}(x,t)+ \dot{s}(t) \phi'(x-s(t))X(t) ,\\
\label{tarBC2} w_x(0,t) =&  \fr{\beta}{\alpha} \phi(0) u(0), \\
\label{tarBC1} w(s(t),t) =& \ep X(t), \\
\label{tarODE}\dot{X}(t)=&-cX(t)-\beta w_x(s(t),t).
\end{align}
The objective of the transformation \eqref{eq:DBST} is to add a stabilizing term $-c X(t)$ in \eqref{tarODE} of the target $(w,X)$-system which is easier to prove the stability than $(u,X)$-system.
By taking the derivative of \eqref {eq:DBST} with respect to $t$ and $x$ respectively, to satisfy \eqref{tarPDE}, \eqref{tarBC1}, \eqref{tarODE}, we derive the conditions on the gain kernel solution, and they leads to the following solution:
\begin{align}
\label{kernel}\phi(x) =& \frac{c}{\beta} x- \ep.
\end{align}
By taking the derivative of the transformation \eqref{eq:DBST} in $x$ and substituting $x = 0$,  we have
\begin{align}
w_{x}(0,t)=&- \frac{q_{\rm c}(t)}{k} -  \frac{\beta}{\alpha} \ep u(0,t) -\frac{c}{\alpha} \int_{0}^{s(t)} u(y,t) {\rm d}y- \frac{c}{\beta} X(t). \label{tarBC_middle} 
\end{align}
Substituting the design of the sampled-data control $q_{\rm c}(t) = q_{j} = - c \widetilde E_{j}$ for all $ t\in[t_{j}, t_{j+1})$ and for all $j \in {\mathcal Z}^+$, and recalling the definition of $\widetilde E(t)$ in \eqref{Energy}, the boundary condition \eqref{tarBC_middle} can be written as
\begin{align}\label{tarBC_middle2}
w_{x}(0,t)=& -\frac{c}{k} \left( \widetilde E(t) - \widetilde E_{j}\right) -  \frac{\beta}{\alpha} \ep u(0,t) . 
\end{align}
Moreover, substituting \eqref{EtEk}, we can describe \eqref{tarBC_middle2} as
\begin{align}
w_{x}(0,t)=& f(t) -  \frac{\beta}{\alpha} \ep u(0,t) , \label{wx0-semi} 
\end{align}
where $f(t)$ is an explicit function in time defined by 
\begin{align}
f(t) = \frac{c^2}{k} \widetilde 
E_{j} \cdot (t - t_{j} )  , \quad \forall t\in[t_{j}, t_{j+1}), \quad j \in {\mathcal Z}^+. 
\end{align}
The closed form representation of \eqref{wx0-semi} using variables $(w,X)$ is given after the inverse transformation is obtained in the next section. 

\subsubsection{Inverse transformation}
Consider the following inverse transformation
\begin{align}\label{inv-trans}
u(x,t)=&w(x,t)-\frac{\beta}{\alpha} \int_{x}^{s(t)} \psi (x-y)w(y,t) {\rm d}y \notag\\
&-\psi(x-s(t)) X(t).
\end{align}
Taking the derivatives of \eqref{inv-trans} in $x$ and $t$ along \eqref{tarPDE}-\eqref{tarODE}, we obtain the gain kernel solution as
\begin{align}\label{inv-gain}
\psi(x) =& e^{ \lambda x } \left( p \sin\left( \omega x \right) + \ep \cos\left( \omega x \right) \right) ,
\end{align}
where $\lambda = \fr{\beta \varepsilon}{2 \alpha}$, $ \omega = \sqrt{\fr{4 \alpha c - (\varepsilon\beta)^2 }{4 \alpha^2 } }$, $p = - \fr{1}{2 \alpha \beta \omega} \left( 2 \alpha c - (\varepsilon \beta )^2 \right) $,
and $0<\ep<2 \fr{\sqrt{\alpha c}}{\beta}$ is to be chosen later. Finally, using the inverse transformation, the boundary condition \eqref{tarBC2} is rewritten as
\begin{align}
\label{eq:tarBC2} w_x(0,t) =&f(t)  - \frac{\beta}{\alpha}\varepsilon \left[  w(0,t) \right. \notag\\
& \left. -\frac{\beta}{\alpha} \int_{0}^{s(t)} \psi (-y)w(y,t) {\rm d}y-\psi(-s(t)) X(t) \right].
\end{align}
Therefore, the closed form of the target $(w,X)$-system is described by \eqref{tarPDE}, \eqref{tarBC1}, \eqref{tarODE}, and \eqref{eq:tarBC2}. 
%
%
\subsubsection{Lyapunov method}
To show the stability of the original system, first we show the stability of the target system \eqref{tarPDE}, \eqref{tarBC1}, \eqref{tarODE}, and \eqref{eq:tarBC2}. For a given $t \geq 0$, we define the most recent sampling number as 
\begin{align} 
n:= \{n \in {\mathcal Z}^+| t_{n} \leq t < t_{n+1}\} , 	
\end{align}
and we firstly apply Lyapunov method for the time interval $t \in [t_{j}, t_{j+1})$ for all $j =0,1, \cdots, n-1$, and next for the interval from $t_{n}$ to $t$. For both cases, we consider the following functional
\begin{align}\label{lyap}
V = \fr{1}{2\alpha } || w||^2 + \fr{\varepsilon}{2\beta } X(t)^2,
\end{align}
where $||w||$ denotes $L_2$ norm defined by $||w|| = \sqrt{ \int_0^{s(t)} w(x,t)^2 {\rm d}x } $. Note that Poincare's and Agmon's inequalities for the system \eqref{tarPDE}--\eqref{tarBC1} with $0<s(t)<s_r$ lead to 
\begin{align} \label{Poincare}
	|| w||^2 \leq 2 s_r \ep^2 X(t)^2 + 4 s_r^2 || w_{x}||^2, \\
	w(0,t)^2 \leq 2 \ep^2 X(t)^2 + 4 s_r  || w_{x}||^2. \label{Agmon} 
\end{align}
Taking the time derivative of \eqref{lyap} along with the solution of   \eqref{tarPDE}--\eqref{tarODE}, \eqref{eq:tarBC2}, we have
\begin{align}\label{Vdot}
\dot{V} =&  - || w_{x}||^2  - \fr{\varepsilon}{\beta}cX(t)^2 - w(0,t) f(t) +\frac{\beta}{\alpha}\varepsilon w(0,t)^2 \notag\\
& -\frac{\beta}{\alpha}\varepsilon w(0,t)\left[ \frac{\beta}{\alpha} \int_{0}^{s(t)} \psi (-y)w(y,t) {\rm d}y +\psi(-s(t)) X(t) \right] \notag\\
&+ \fr{\dot{s}(t)}{\alpha} \left( \fr{\ep^2}{2 } X(t)^2+  \fr{ c}{\beta} \int_0^{s(t)} w(x,t) {\rm d}x X(t)  \right) .
\end{align}
Applying Young's inequality to the second line of \eqref{Vdot} twice, we get
 \begin{align}\label{dyoung}
& - w(0,t) f(t) \leq \frac{\gamma_1}{2} w(0,t)^2 + \frac{1}{2\gamma_1} f(t)^2 , \\
 & - w(0,t)\left[\frac{\beta}{\alpha} \int_{0}^{s(t)} \psi (-y)w(y,t) {\rm d}y+\psi(-s(t)) X(t) \right] \notag\\
  \leq& \fr{1}{2}w(0,t)^2 + \frac{\beta^2}{\alpha^2 \gamma_2} \left(\int_{0}^{s(t)} \psi (-y)w(y,t) {\rm d}y \right)^2 \notag\\
  &+ \gamma_2 \psi(-s(t))^2 X(t)^2, \label{yyoung2}
  \end{align}
 where $\gm_{1}>0$ and $\gm_2>0$ are parameters to be determined. Applying \eqref{dyoung}, \eqref{yyoung2}, \eqref{Poincare}, \eqref{Agmon}, and Cauchy Schwarz inequalities to \eqref{Vdot} with choosing $\gm_1 = \frac{1}{4 s_r} $ and $\gm_2 = \fr{1}{8}$, we have
\begin{align}\label{Vdot3}
\dot{V} \leq &  - \left(\fr{1}{2} - \fr{2 \beta s_r}{\alpha} \left( \frac{ 64 c s_r^2}{\alpha } + 3 \right)\varepsilon\right) || w_{x}||^2 \notag\\
& - \varepsilon \left(\fr{c}{8\beta} + g(\ep) \right)X(t)^2 + 2 s_r f(t)^2    \notag\\
&+ \fr{\dot{s}(t)}{2\alpha} \left( \ep^2 X(t)^2 +  \fr{ 2c}{\beta} \left|\int_0^{s(t)} w(x,t) {\rm d}x X(t)  \right|\right) ,
\end{align}
where $g(\ep) = \fr{c}{8\beta} - \fr{\ep}{4 s_r} -  \fr{ \beta }{\alpha} \left( \frac{ 64 c s_r^2}{\alpha } + 3 \right) \ep^2 $. Since $g(0) = \fr{c}{8\beta }>0$ and $g'(\ep) = - \fr{1}{4 s_r} - \frac{2\beta \ep }{\alpha}  \left( \frac{ 64 c s_r^2}{\alpha } + 3 \right) <0$ for all $ \ep>0$, there exists $\ep^*$ such that $g(\ep)>0$ for $0<\ep<\ep^*$ and $g(\ep^*)=0$. Thus, setting $\ep < \min\left\{ \ep^*, \fr{\alpha}{8 \beta s_r \left( \frac{ 64 c s_r^2}{\alpha } + 3 \right)} \right\}$ ,
the inequality \eqref{Vdot3} leads to
\begin{align}\label{dotV}
\dot{V} \leq &  - b V + 2 s_r f(t)^2  + a\dot{s}(t) V , 
\end{align}
where 
\begin{align} \label{ab-def} 
b =\fr{1}{8} \min \left\{ \fr{\alpha}{s_r^2}, c \right\}, \quad a = \fr{2\beta \ep }{\alpha} \max \left\{1,\fr{\alpha c^2 s_r}{2\beta^3 \ep^3} \right\}.
\end{align} 
Consider the following functional
\begin{align}\label{Wdef}
W = V e^{ - a s(t)} .
\end{align}
Taking the time derivative of \eqref{Wdef} with the help of \eqref{dotV}, we deduce
\begin{align}
\dot{W} & \leq - b W + 2 s_r f(t)^2 e^{-a s(t)} \notag\\
&\leq - b W + 2s_r f(t)^2. \label{Wdot}
\end{align}

\textbf{(i) For $t \in [t_{j}, t_{j+1})$, for all $j=0,1, \cdots, n-1,$}\\
Applying comparison principle to \eqref{Wdot} for $t \in [t_{j}, t_{j+1})$ leads to
\begin{align}\label{Wineq}
W(t) \leq W(t_{j}) e^{-b(t-t_j)} + 2 s_r e^{-b t} \int_{t_{j}}^{t} e^{b \tau}f(\tau)^2 {\rm d}\tau .
\end{align}
Setting $t = t_{j+1}$ and recalling $f(t) = \frac{c^2}{k} \widetilde E_{j}  (t - t_{j} ), \forall t\in[t_{j}, t_{j+1})$, we get 

\begin{align} 
	W_{j+1} \leq & W_{j} e^{-b \tau_{j}} + \frac{2 c^4 s_r}{k^2} e^{-b \tau_{j}} \widetilde E_{j}^2 I_j ,  \label{Wjbound1} 
	\end{align} 
where $W_{j} = W(t_j)$, and $I_j$ is defined by 
\begin{align} 
I_j := \int_{t_{j}}^{t_{j+1}} e^{b (\tau-t_{j})}(\tau-t_j)^2 {\rm d}\tau . 
\end{align} 
Then, by introducing the variable $s = b (\tau - t_{j})$ and integration by substitution, with the help of $b \tau_j < \fr{1}{8} c \tau_j < \fr{1}{8}$ for all $j \in {\mathcal Z}^+$ derived by \eqref{ab-def}, Assumption \ref{ass:sample} and the fact that $c<\fr{1}{R}$, one can derive the following inequality: 
\begin{align} 
I_j = & \fr{1}{b^3} \int_{0}^{b \tau_j } e^{s}s^2 {\rm d}s \leq  \fr{J}{b^3}, \label{Ij-J} 
\end{align} 
where $J$ is defined by $J:= \int_{0}^{\fr{1}{8} } e^{s}s^2 {\rm d}s $. Applying \eqref{Ij-J} to \eqref{Wjbound1} yields 
\begin{align} 
W_{j+1}	\leq & W_{j} e^{-b \tau_{j}} + B_{j}, \label{Wj+1}
\end{align} 
where $B_j$ is defined by 
\begin{align} 
B_j = 	\frac{2 Jc^4 s_r}{k^2 b^3} e^{-b \tau_{j}} \widetilde E_{j}^2. \label{Bj-def} 
\end{align}
Applying \eqref{Wj+1} from $j=n-1$ to $j=0$ inductively, we get 
\begin{align} 
W_{n} \leq 	W_0 e^{- b \sum_{i=0}^{n-1} \tau_{i}} + B_{n-1} + \sum_{i=0}^{n-2} B_{i} e^{-b \sum_{j=i+1}^{n-1} \tau_{j}} .  \label{Wnbound-1} 
\end{align}
By \eqref{Bj-def} and the solution of $\widetilde E_{j}$ given in \eqref{Ejsol}, we have 
\begin{align}
&\sum_{i=0}^{n-2} B_{i} e^{-b \sum_{j=i+1}^{n-1} \tau_{j}} \notag\\
\leq &	\frac{2 Jc^4 s_r \widetilde E_0^2 e^{-b \sum_{j=0}^{n-1} \tau_{j}}}{k^2 b^3} \left( 1+  \sum_{i=1}^{n-2} \left(\prod_{k=0}^{i-1} \left( 1 - c \tau_k \right)^2 \right) e^{b \sum_{j=0}^{i-1} \tau_{j}} \right) \notag\\
\leq &	\frac{2 Jc^4 s_r \widetilde E_0^2 e^{-b \sum_{j=0}^{n-1} \tau_{j}}}{k^2 b^3} \left( 1+  \sum_{i=1}^{n-2} \left(\prod_{k=0}^{i-1} \left( 1 - c \tau_k \right)^2 e^{b \tau_{k}}\right)  \right) . \label{delta-ineq}
\end{align}
Since $b =  \fr{1}{8} \min \left\{ \fr{\alpha}{s_r^2}, c \right\} <  \frac{c}{8}$, by using $r = \inf_{j \in {\mathcal Z}^+} \{\tau_j\} > 0$ given in Assumption \ref{ass:sample}, the following inequality holds 
\begin{align} 
\left( 1 - c \tau_{i} \right)^2   e^{b \tau_{i}} \leq \left( 1 - c r \right)^2   e^{\frac{ cr}{8}}  := \delta < 1, \quad \forall j \in {\mathcal Z}^+.  
\end{align} 
Thus, the inequality \eqref{delta-ineq} leads to 
\begin{align}
\sum_{i=0}^{n-2} B_{i} e^{-b \sum_{j=i+1}^{n-1} \tau_{j}} \leq &	\frac{2 Jc^4 s_r \widetilde E_0^2 e^{-b \sum_{j=0}^{n-1} \tau_{j}}}{k^2 b^3} \left( 1+  \sum_{i=1}^{n-2} \delta^i  \right) \notag\\
\leq & \frac{2 Jc^4 s_r \widetilde E_0^2 }{k^2b^3 (1 - \delta)} e^{-b \sum_{j=0}^{n-1} \tau_{j}} .  \label{Bi-bound}
\end{align}
In the similar way, we get 
\begin{align} \label{Bn-1sol} 
B_{n-1} \leq 	\frac{2 Jc^4 s_r \widetilde E_0^2}{k^2 b^3 (1 - \delta)} e^{-b \sum_{j=0}^{n-1} \tau_{j}} . 
\end{align}
Recalling that $\tau_j = t_{j+1} - t_{j}$ and $t_0 = 0$, we get $\sum_{j=0}^{n-1} \tau_{j} = t_n$. Applying \eqref{Bi-bound} and \eqref{Bn-1sol} to \eqref{Wnbound-1}, we arrive at 
\begin{align} 
	W_{n} \leq 	(W_0  + A \widetilde E_0^2  )e^{-b t_{n}} . \label{Wn-bound} 
\end{align}
 where $A = \frac{2 Jc^4 s_r}{k^2 b^3 (1 - \delta)}$. 
 
 \textbf{(ii) For $t \in [t_{n}, t_{n+1})$,}\\
 Applying comparison principle to \eqref{Wdot} from $t_{n}$ to $t \in [t_{n}, t_{n+1})$, we get 
 \begin{align} 
 W(t) \leq & W_{n} e^{-b(t-t_n)} + B_{n} e^{-b(t - t_{n+1})} \notag\\
 \leq & 	 W_{n} e^{-b(t-t_n)} + A \widetilde E_0^2 e^{-bt} . \label{Wt-bound} 
 \end{align}
Finally, combining \eqref{Wn-bound} and \eqref{Wt-bound}, the following bound is obtained 
\begin{align} 
W(t) \leq (W_0 + 2 A \widetilde E_0^2) 	e^{-bt} . \label{Wt-exp} 
\end{align}
Recalling the relation $W = V e^{-a s(t)}$ defined in \eqref{Wdef}, and applying $0<s(t)<s_r$, the norm estimate for $W$ in \eqref{Wt-exp} leads to the following estimate for $V$: 
\begin{align} \label{Vt-ineq-norm} 
V(t) \leq e^{ a s_r}  (V_0 + 2 A \widetilde E_0^2) e^{-bt} .
\end{align}
We consider the $L_2$-norm of $(u,X)$-system defined by 
\begin{align} 
	\Psi(t) =  \int_{0}^{s(t)} u(x,t)^2 {\rm d}x + X(t)^2. 
\end{align}
Due to the invertibility of the transformation from $(u,X)$ to $(w,X)$ together with the boundedness of the domain $0<s(t)<s_r$, there exist positive constants $\underline{M}>0$ and $\overline{M}>0$ such that the following inequalities hold: 
\begin{align}\label{Psibound} 
\underline{M} \Psi(t) \leq V(t) \leq \overline{M} \Psi(t). 	
\end{align}
Moreover, due to the definition of the reference energy $\widetilde E(t) = \fr{k}{\alpha } \int_0^{s(t)} u(x,t) {\rm d}x + \fr{k}{\beta} X(t)$ given in \eqref{Energy}, using Young's and Cauchy Schwarz inequalities one can show that 
\begin{align} \label{E0Psi0} 
\widetilde E_0^2 \leq K \Psi_0, 
\end{align}
where $K = 2 k^2 \max\{\frac{s_r}{\alpha^2}, \frac{1}{\beta^2}\}$. Applying \eqref{Psibound} and \eqref{E0Psi0} to \eqref{Vt-ineq-norm}, we deduce that there exists positive constant $M>0$ such that the following inequality holds 
\begin{align} 
\Psi(t) \leq M \Psi_0 e^{-bt}, 	
\end{align}
which completes the proof of Theorem \ref{theo-1}. 

\section{Sampled-Data Design for Two-Phase Stefan Problem} \label{sec:twophase} 
\begin{figure}[t]
\centering
\includegraphics[width=2.5in]{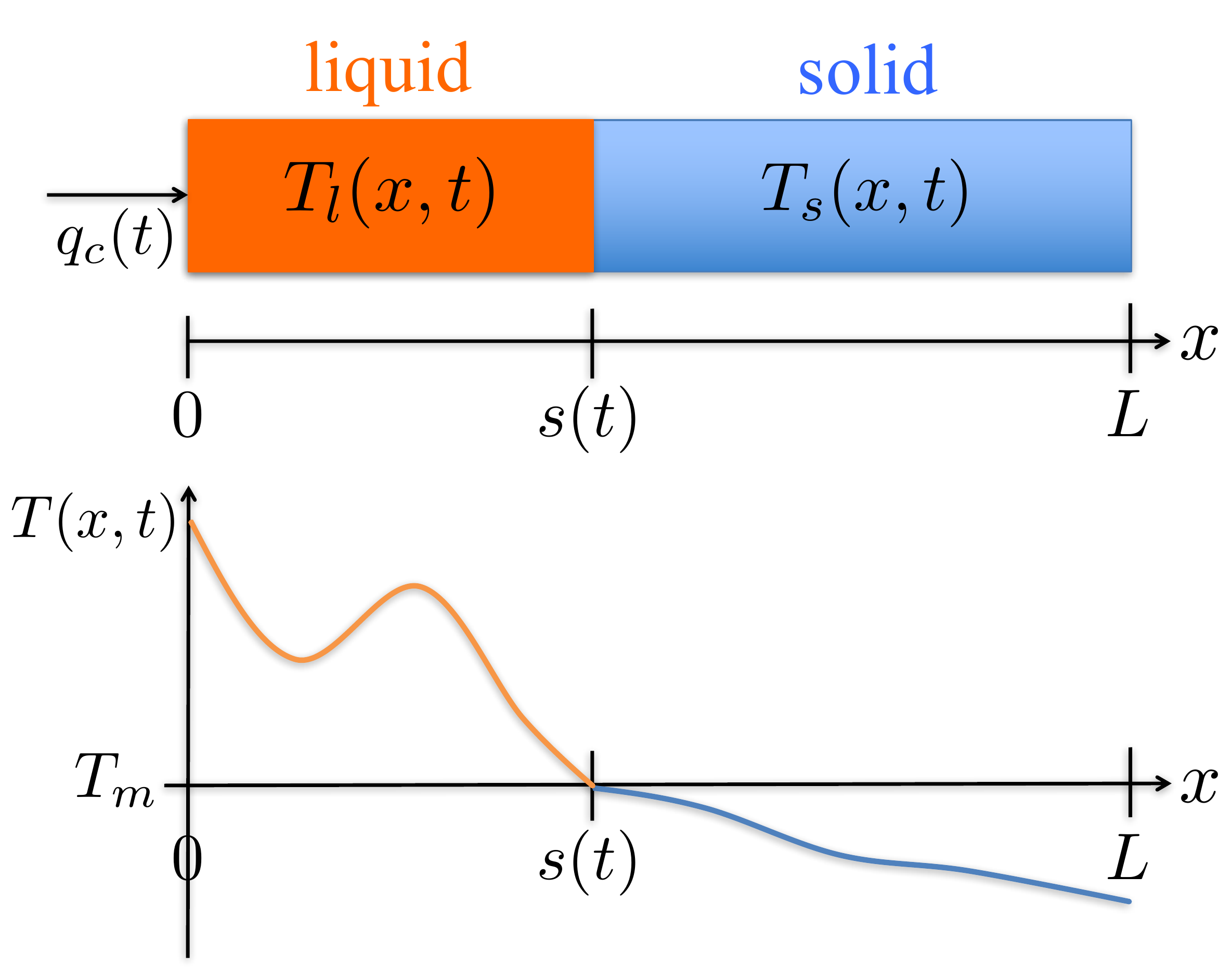}\\
\caption{Schematic of the two-phase Stefan problem.}
\label{fig:2phstefan}
\end{figure}
In this section, we extend the results we have established in the previous section to the "two-phase" Stefan problem, where the temperature dynamics in the solid phase is governed by the heat equation with different physical parameters from the liquid phase, following the work in \cite{koga2018CDC}. This configuration is depicted in Fig. \ref{fig:2phstefan}. 
\subsection{Problem statement} 
The governing equations are descried by the following coupled PDE-ODE-PDE system: 
\begin{align}\label{twoph-sys1}
 \fr{\pa T_{{\rm l}}}{\pa t}(x,t) =&\alpha_{{\rm l}}  \fr{\pa^2 T_{{\rm l}}}{\pa x^2}(x,t), \quad 0<x<s(t),\\
 \fr{\pa T_{{\rm l}}}{\pa x}(0,t) =& -q_{{\rm c}}(t)/k_{{\rm l}}, \quad T_{{\rm l}}(s(t),t) =T_{{\rm m}},\\
\label{twoph-sys3} \fr{\pa T_{{\rm s}}}{\pa t}(x,t) =&\alpha_{{\rm s}}  \fr{\pa^2 T_{{\rm s}}}{\pa x^2}(x,t), \quad s(t)<x<L, \\
\label{twoph-sys4} \fr{\pa T_{{\rm s}}}{\pa x}(L,t) =& 0,\quad T_{{\rm s}}(s(t),t) =T_{{\rm m}}, \\
\label{twoph-sys5} \gm \dot{s}(t) =& - k_{{\rm l}} \fr{\pa T_{{\rm l}}}{\pa x}(s(t),t)+k_{{\rm s}} \fr{\pa T_{{\rm s}}}{\pa x}(s(t),t),
\end{align}
where $\gamma = \rho_{{\rm l}} \Delta H^*$, and all the variables denote the same physical value with the subscript "l" for the liquid phase and "s" for the solid phase, respectively. The solid phase temperature must be lower than the melting temperature, which serves as one of the conditions for the model validity, as stated in the following remark. 
\begin{remark}
To keep the physical state of each phase meaningful, the following conditions must be maintained: 
\begin{align}\label{valid1-2ph}
T_{{\rm l}}(x,t) \geq& T_{{\rm m}}, \quad \forall x\in(0,s(t)), \quad \forall t>0, \\
\label{valid2-2ph}T_{{\rm s}}(x,t) \leq& T_{{\rm m}}, \quad \forall x\in(s(t),L), \quad \forall t>0, \\
\label{valid3-2ph} 0< &s(t)<L, \quad \forall t>0. 
\end{align}
\end{remark}
\begin{lemma} \label{lem:twoph-flux} 
If the solution to \eqref{twoph-sys1}--\eqref{twoph-sys5} satisfies the conditions \eqref{valid1-2ph}--\eqref{valid3-2ph}, then the following properties hold: 
\begin{align}
\fr{\pa T_{{\rm l}}}{\pa x}(s(t),t) \leq 0, \quad \fr{\pa T_{{\rm s}}}{\pa x}(s(t),t) \leq 0, \quad  \forall t \geq 0. 
\end{align} 
\end{lemma} 
The following assumption on the initial data $(T_{{\rm l},0}(x), T_{{\rm s},0}(x), s_0) := (T_{{\rm l}}(x,0), T_{{\rm s}}(x,0), s(0))$ is imposed.  
\begin{assum}\label{initial-2ph} 
$0<s_0<L$,  $T_{{\rm l},0}(x) \geq T_{{\rm m}}$ for all $x \in [0,s_0]$, $T_{{\rm s},0}(x) \leq T_{{\rm m}}$ for all $x \in [s_0, L]$, and $T_{{\rm l},0}(x) $ and $T_{{\rm s},0}(x) $  are continuously differentiable in $x\in[0,s_0]$ and $x\in[s_0, L]$, respectively. 
 \end{assum}
 The following lemma is provided to ensure the conditions of the model validity. 
\begin{lemma}\label{valid} 
Under Assumption \ref{initial-2ph}, and provided that $q_{{\rm c}}(t)$ is a piecewise continuous function with satisfying 
\begin{align} \label{qmqf} 
q_{{\rm c}}(t)\geq 0,\quad  \forall t \in [0, t^*), 
\end{align} 
there exists a finite time $\overline t:= \sup_{t \in (0, t^*)}\{t | s(t) \in (0,L)\}>0$ such that the solution to \eqref{twoph-sys1}--\eqref{twoph-sys5} exists and unique and satisfies the model validity conditions \eqref{valid1-2ph}--\eqref{valid3-2ph} for all $t \in (0,\overline t)$. Moreover, if $t^* = \infty$ and it holds 
\begin{align} \label{validcondition} 
0< \gm s_{\infty} + \int_0^{t} q_{{\rm c}}(s) {\rm d}s < \gm L ,
\end{align}
for all $t \geq 0$, where 
\begin{align}  \label{sinf}
s_{\infty}:= &s_0 + \frac{k_{{\rm l}}}{\alpha_{{\rm l}} \gm} \int_{0}^{s_0} (T_{{\rm l},0}(x) - T_{{\rm m}}) {\rm d}x  \notag\\
&+ \frac{k_{{\rm s}}}{\alpha_{{\rm s}}\gm} \int_{s_0}^{L} (T_{{\rm s},0}(x) - T_{{\rm m}}) {\rm d}x	, 
\end{align}
then $\overline t = \infty$, namely, the well-posedness and the model validity conditions are satisfied for all $t \geq 0$.  
 \end{lemma}
 
Lemma \ref{valid} is proven in \cite{Cannon71flux} (Theorem 1 in p.4 and Theorem 4 in p.8) by employing the maximum principle. The variable $s_{\infty}$ defined in \eqref{sinf} is the final interface position $s_{\infty} = \lim_{t \to \infty} s(t)$ under the zero input $q_{{\rm c}}(t) \equiv 0$ for all $t \geq 0$.  For \eqref{validcondition} to hold for all $t \geq 0$, we at least require it to hold at $t = 0$, which leads to the following assumption. 
\begin{assum}\label{ass:E0}
The variable $s_{\infty}$ defined in \eqref{sinf} given by initial values satisfies
\begin{align} 
0< s_{\infty} < L . 
\end{align}
\end{assum}

\subsection{Control design and main result} 
We apply ZOH to the boundary control design for the two-phase Stefan problem developed in \cite{koga2018CDC}, resulting in the following sampled-data control   
\begin{align}\label{control-2ph}
q_{ c}(t)=& -c \left( \frac{k_{{\rm l}}}{\alpha_{{\rm l}}} \int_{0}^{s(t_{j})} (T_{{\rm l}}(x,t_{j}) - T_{{\rm m}}) {\rm d}x \right. \notag\\
&\left. + \frac{k_{{\rm s}}}{\alpha_{{\rm s}}} \int_{s(t_{j})}^{L} (T_{{\rm s}}(x,t_{j}) - T_{{\rm m}}) {\rm d}x + \gamma (s(t_{j}) - s_{{\rm r}}) \right), 
\end{align}
where $c>0$ is the controller gain, for all $t \in [t_{j}, t_{j+1})$ for all $j \in {\mathcal Z}^+$. The restriction on the setpoint $s_r$ for the two-phase Stefan problem is given by the following. 
\begin{assum}\label{twoph-assum2}
The setpoint is chosen to satisfy
\begin{align} \label{setpoint}
s_{\infty} < s_{{\rm r}} < L,  
\end{align} 
where $s_{\infty}$ is defined in \eqref{sinf}. 
\end{assum}

We state the following theorem for the sampled-data control of the two-phase Stefan problem.

\begin{thm}\label{theo-2}
Consider the closed-loop system \eqref{twoph-sys1}--\eqref{twoph-sys5} and the sampled-data control law \eqref{control-2ph} under Assumptions \ref{initial-2ph}--\ref{twoph-assum2}. Then for every $0<r \leq R<1/c$, there exists a constant $M:=M(r)$ for which the following property holds: for every sequence $\{t_j\geq 0 :  j = 0, 1, 2, \dots \}$ with $t_0=0$ for which Assumption 3 holds, the initial-boundary value problem \eqref{twoph-sys1}--\eqref{twoph-sys5} with \eqref{control-2ph} has a unique solution satisfying \eqref{valid1-2ph}--\eqref{valid3-2ph} as well as the following estimate:
\begin{align}\label{twoph-h1}
\Psi(t) \leq M \Psi(0) \exp{(- \bar bt)} ,
\end{align}
where $\bar b=\fr{1}{8} \min \left\{ \fr{\alpha_{\rm l}}{L^2}, \fr{4\alpha_{\rm s}}{L^2}, c \right\}$, for all $t\geq0$, in the $L_2$ norm $\Psi(t) =  \int_{0}^{s(t)} \left(T_{\rm l}(x,t)-T_{\rm m}\right)^2 {\rm d}x +  \int_{s(t)}^{L} \left(T_{\rm s}(x,t)-T_{\rm m}\right)^2 {\rm d}x +(s(t)-s_{r})^2$.

\end{thm}

As in the previous section, the equivalence of the closed-loop system under the control law \eqref{control-2ph} with the system under an open-loop input is presented in the following lemma. 
\begin{lem} \label{lem:2ph} 
The closed-loop system consisting of \eqref{twoph-sys1}--\eqref{twoph-sys5} with the control law \eqref{control-2ph} has a unique classical solution satisfying \eqref{valid1-2ph}--\eqref{valid3-2ph}, which is equivalent to the open-loop solution of \eqref{twoph-sys1}--\eqref{twoph-sys5} with 
\begin{align} \label{qjsol-2ph}
	q_{\rm c}(t) = q_{j} = & q_{0}\prod_{i=0}^{j-1} \left( 1 - c \tau_i \right), \quad \forall t \in [t_{j}, t_{j+1}), \quad \forall j \in {\mathcal Z}^{+}, 
\end{align}
where 
\begin{align} \label{qini-2ph}
q_0 = - c &\bigg(\frac{k_{\rm l}}{\alpha_{\rm l}} \int_0^{s_0} (T_{\rm l,0}(x) - T_{{\rm m}}) {\rm d}x \notag\\
& + \frac{k_{\rm s}}{\alpha_{\rm s}} \int_0^{s_0} (T_{\rm s,0}(x) - T_{{\rm m}}) {\rm d}x + \gm (s_0 - s_r ) \bigg).    
\end{align}
\end{lem} 
\begin{pf} 
The proof of Lemma \ref{lem:2ph} is almost same procedure as the proof of Lemma \ref{lem2} once we redefine the system's internal energy as 
\begin{align}\label{Et-2ph}
\widetilde E(t) =& \frac{k_{{\rm l}}}{\alpha_{{\rm l}}} \int_{0}^{s(t)} (T_{{\rm l}}(x,t) - T_{{\rm m}}) {\rm d}x \notag\\
&+ \frac{k_{{\rm s}}}{\alpha_{{\rm s}}} \int_{s(t)}^{L} (T_{{\rm s}}(x,t) - T_{{\rm m}}) {\rm d}x + \gm s(t) , 
\end{align}
and obtain the same differential equation of the energy as \eqref{conservation}. Since the control law \eqref{control-2ph} is equivalent to $q_{\rm c}(t) = q_{j} = - c \widetilde E_{j}$ for all $t \in [t_{j}, t_{j+1})$ and for all $j \in {\mathcal Z}^+$, in the same manner as the proof of Lemma \ref{lem2}, one can derive the equivalence of the solution to the open-loop solution with \eqref{qjsol-2ph} and \eqref{qini-2ph}. By applying Lemma \ref{lem:2ph}, the well-posedness of the solution is proven with satisfying the conditions \eqref{valid1-2ph}--\eqref{valid3-2ph} for the model validity.  
\end{pf}

To prove the exponential stability estimate \eqref{twoph-h1}, by following the procedure in \cite{koga2018CDC}, first we introduce the reference error states as follows. 
\begin{align}\label{ref-2ph}
u(x,t) :=& T_{\rm l} (x,t)-T_{{\rm m}}, \\
v(x,t) : = & T_{\rm s} (x,t)-T_{{\rm m}}, \\
X(t) :=& s(t) - s_{{\rm r}} + \fr{\beta_{{\rm s}}}{\alpha_{{\rm s}}} \int_{s(t)}^{L} v(x,t) {\rm d}x . \label{refX-2ph} 
\end{align}
Using these reference error variables, the total PDE-ODE-PDE system given in \eqref{twoph-sys1}--\eqref{twoph-sys5} is reduced to the following PDE-ODE system
\begin{align}\label{uX-sys1}
 u_{t}(x,t) =& \alpha_{{\rm l}} u_{xx}(x,t), \quad 0<x<s(t), \\
\label{uX-BC}u_x(0,t) = &- q_{{\rm c}}(t)/k_{{\rm l}}, \quad  u(s(t),t) = 0,\\
\dot{X}(t) =& - \beta_{{\rm l}} u_x(s(t),t) .  \label{uX-sys3}
\end{align}
Note that the formulation of the above system is equivalent to \eqref{u-sys1}--\eqref{u-sys4} which is the reference error system in the one-phase case. The only difference is the non-monotonic property of the moving interface, namely, $\dot s(t) \geq 0$ is no longer verified, which is utilized for the stability proof in Section \ref{sec:stability}. However, we can deal with the problem by following the procedure in \cite{koga2018CDC}. Owing to the properties $u_{x}(s(t),t) <0$ and $v_x(s(t),t)<0$ derived from Lemma \ref{lem:twoph-flux}, it holds that $|\dot s(t) | \leq - \beta_{\rm l} u_{x}(s(t),t) - \beta_{\rm s} v_{x}(s(t),t)$. Introduce 
\begin{align} \label{ztdef}
z(t) := X(t)+ \fr{\beta_{\rm s}}{\alpha_{{\rm s}}} \int_{s(t)}^{L} v(x,t) {\rm d}x . 
\end{align}  
Then, due to the negativities $X(t) <0$ as derived in \eqref{Xtnegative} and $v(x,t) < 0$ deduced from \eqref{valid2-2ph}, it holds $z(t)<0$. Moreover, taking the time derivative of \eqref{ztdef} leads to 
\begin{align} 
\dot{z}(t) = - \beta_{\rm l} u_{x}(s(t),t) - \beta_{\rm s} v_x(s(t),t) >0. 
\end{align} 
Hence, owing to Assumption \ref{initial-2ph}, there exists a positive constant $\delta >0$ such that $- \delta < z(0) < z(t)<0$ holds. Therefore, following the same procedure as in Section \ref{sec:stability} with replacing the term $\dot s(t)$ with $\dot z(t)$, it is straightforward to derive that there exists a positive constant $N>0$ such that the following norm estimate holds: 
\begin{align} \label{phiineq-2ph}
\Phi(t) \leq & N \Phi(0) e^{- b t}, 
\end{align} 
where $\Phi(t) = \left( \int_{0}^{s(t)} u(x,t)^2 {\rm d}x + X(t)^2 \right)^{\frac{1}{2}}$, $ b =\fr{1}{8} \min \left\{ \fr{\alpha_{{\rm l}}}{L^2}, c \right\} $. 
Let us define the following three functionals 
\begin{align}\label{V1def-2ph}
V_1 (t)=& \int_0^{s(t)}  (T_{{\rm l}}(x,t) - T_{{\rm m}})^2 {\rm d}x, \\
V_2 (t)=& \int_{s(t)}^{L} (T_{{\rm s}}(x,t) - T_{{\rm m}})^2 {\rm d}x, \label{V2def-2ph} \\
V_3 (t)=& (s(t) - s_{{\rm r}})^2. \label{V3def-2ph}
\end{align} 
Taking the time derivative of \eqref{V2def-2ph} along with the solid phase dynamics \eqref{twoph-sys3} and \eqref{twoph-sys4}, we get
\begin{align} 
 \dot{V}_2 (t)= & -2 \alpha_s \int_{s(t)}^{L} \left(\frac{\pa T_{{\rm s}}}{\pa x}(x,t) \right)^2 {\rm d}x  . \label{V2dot} 
\end{align} 
Applying Young's, Cauchy-Schwarz, Poincare's and Agmon's inequalities to \eqref{V2dot}, we arrive at the following differential inequality 
\begin{align} \label{Vdot2ineq-2ph} 
\dot{V}_2 (t)\leq - \frac{\alpha_{{\rm s}}}{2 L^2} V_2(t)  . 
\end{align} 
Applying the comparison principle to \eqref{Vdot2ineq-2ph}, one can derive 
\begin{align} \label{V2ineq-2ph}
V_2 (t) \leq V_2(0) e^{-  \frac{\alpha_{{\rm s}}}{2 L^2} t}  . 
\end{align} 
Taking the square of \eqref{refX-2ph}, and applying Young's and Cauchy-Schwarz inequalities with the help of $0 < s(t) < L$, one can obtain the following inequality, 
\begin{align} \label{Xineq}
 X(t)^2 \leq & 2 V_3(t) + \fr{2 L \beta_{{\rm s}}^2}{\alpha_{{\rm s}}^2} V_2 (t). 
\end{align}
Applying the same manner to the relation $s(t) - s_{{\rm r}} = X(t) - \frac{\beta_s}{\alpha_s} \int_{s(t)}^{L} (T_{{\rm s}}(x,t) - T_{{\rm m}}) {\rm d}x$ obtained by \eqref{refX-2ph}, one can also derive 
\begin{align}\label{Ytineq}
 V_3(t)   \leq  2 X(t)^2 + \fr{2 L \beta_{{\rm s}}^2}{\alpha_{{\rm s}}^2} V_2 (t). 
\end{align}
Combining \eqref{phiineq-2ph}, \eqref{V2ineq-2ph}, \eqref{Xineq}, and \eqref{Ytineq} using the definitions in \eqref{V1def-2ph}--\eqref{V3def-2ph}, the estimate of the norm $\Psi(t) = V_1(t) + V_2(t) + V_3(t)$ is obtained by the inequality \eqref{twoph-h1} for some positive constant $M>0$, which completes the proof of Theorem \ref{theo-2}.

\section{Numerical Simulation} \label{sec:simulation} 

	\begin{table}[t]
	\vspace{2mm}
	
\caption{Physical properties of paraffin (liquid)}
\begin{center}
    \begin{tabular}{| l | l | l | }
    \hline
    $\textbf{Description}$ & $\textbf{Symbol}$ & $\textbf{Value}$ \\ \hline
    Density & $\rho$ & 790 ${\rm kg}\cdot {\rm m}^{-3}$\\ 
    Latent heat of fusion & $\Delta H^*$ & 210 ${\rm J}\cdot {\rm g}^{-1}$ \\ 
    Heat Capacity & $C_{{\rm p}}$ & 2.38 ${\rm J} \cdot {\rm g}^{-1}\cdot$$^\circ$C$^{-1}$  \\ 
    Melting Temperature & $T_{\rm m}$ & 37.0 $^\circ$C  \\  
    Thermal conductivity & $k$ & 0.220 ${\rm W}\cdot {\rm m}^{-1}$  \\ \hline
    \end{tabular}
\end{center}
\end{table}

Simulation results are performed for the one-phase Stefan problem by considering a cylinder of paraffin whose physical parameters are given in Table 1. Here, we use the well known boundary immobilization method combined with finite difference semi-discretization \cite{kutluay97}. The setpoint and the initial values are chosen as $s_{{\mathrm r}}$ = 2.0 cm, $s_0$ = 0.1 cm, and $T_0(x)-T_{{\mathrm m}}= \bar{T}_0(1-x/s_0)$ with $ \bar{T}_0$ = 1 $^\circ$C. Then, the setpoint restriction stated in Assumption \ref{ass1} is satisfied. We consider periodic sampling with period given by
\begin{align} 
\tau_j = R = \textrm{10 [min]}, \quad \forall j \in {\mathcal Z}	. 
\end{align}
The control gain is set as $c =$ 5.0 $\times$ 10$^{-3}$/s, by which the requirement $R < \fr{1}{c}$ is satisfied. 

\begin{figure}[t]
\centering 
\subfloat[Convergence of the interface to the setpoint is observed without the overshoot.]
{\includegraphics[width=2.7in]{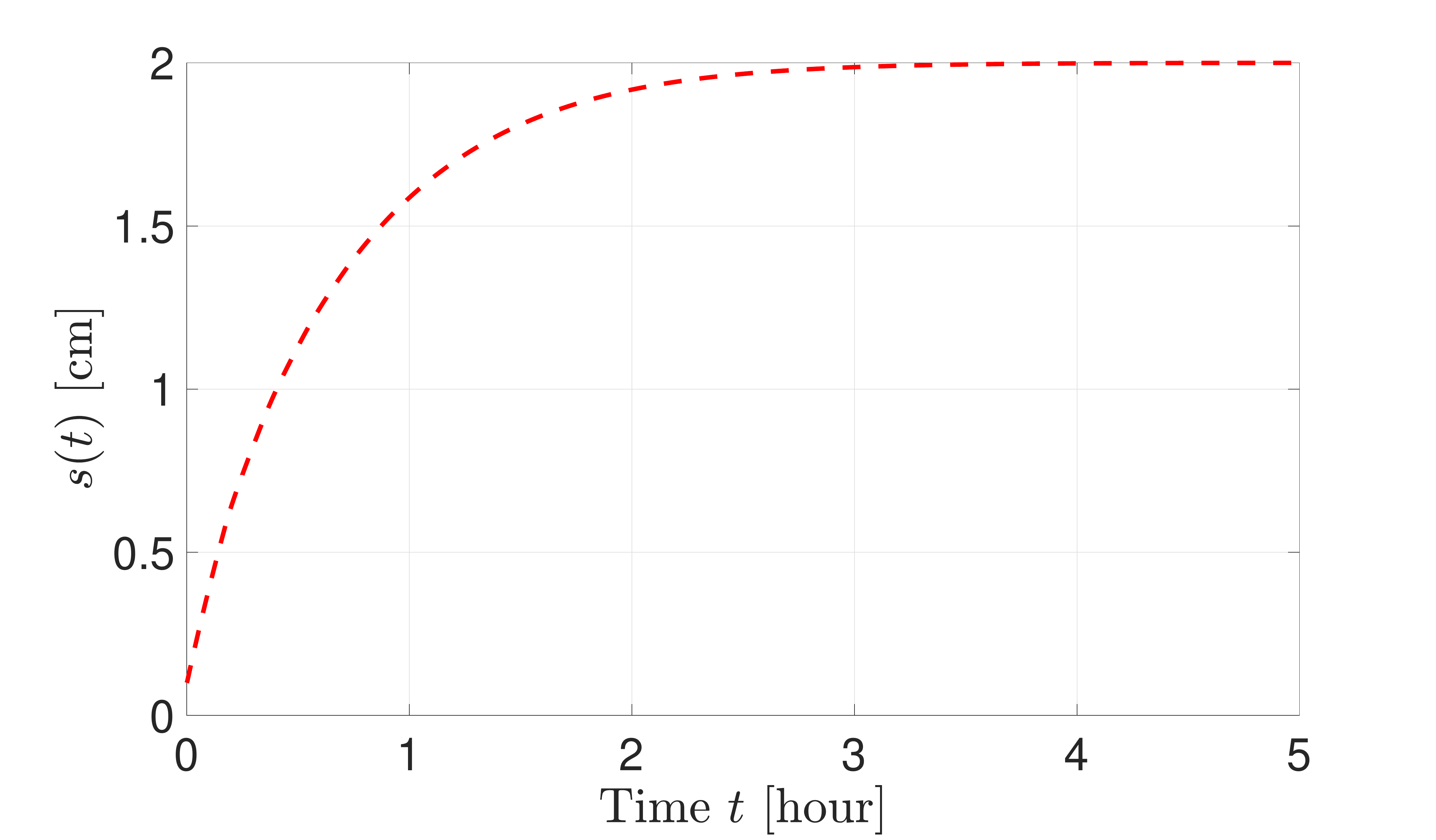}\label{fig:interface}}\\
\subfloat[Positivity of the closed-loop controller is satisfied.]
{\includegraphics[width=2.7in]{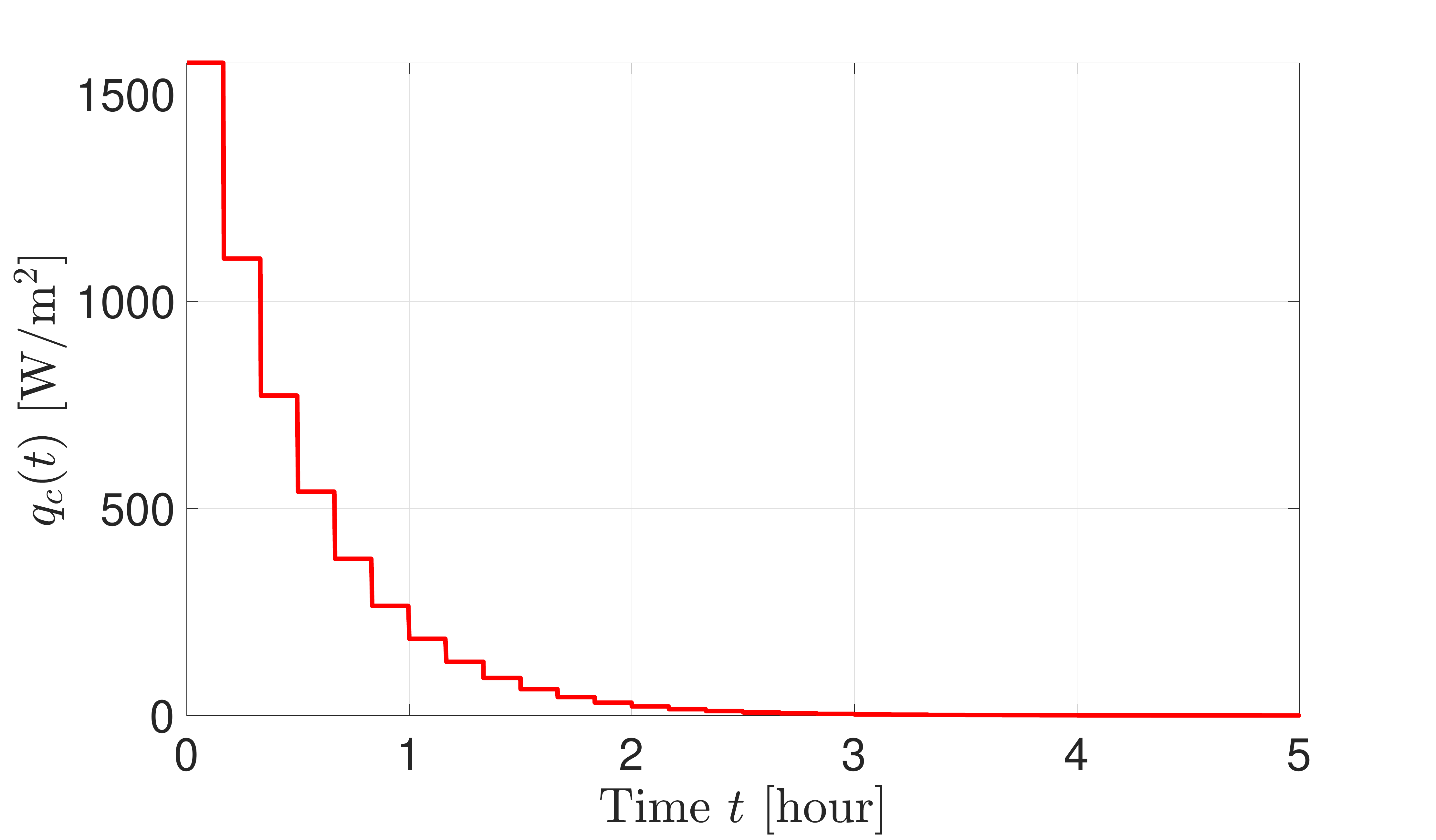}\label{fig:qc}}\\
\subfloat[The model validity of the boundary liquid temperature holds, i.e., $T(0,t)>T_{{\rm m}}$.]
{\includegraphics[width=2.7in]{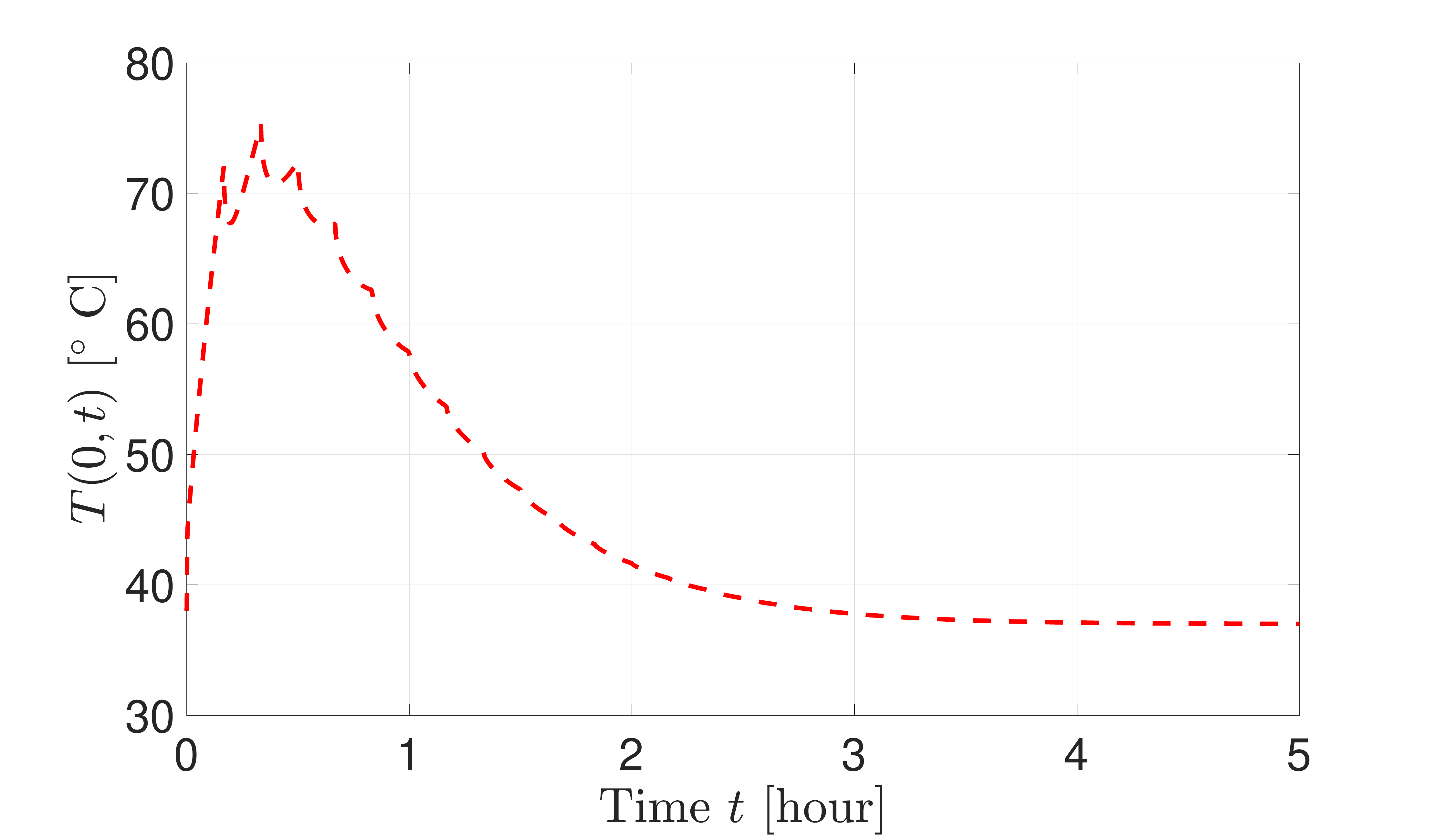}\label{fig:T0}}
\caption{The responses of the system \eqref{sys1}--\eqref{sys4} with ZOH-based sampled-data control \eqref{sampledcont}.  }
\label{fig:response}
\end{figure}

The time responses of the interface position, the control input, and the boundary temperature under the closed-loop system are depicted in Fig. \ref{fig:response} (a)--(c), respectively. Fig. \ref{fig:response} (a) illustrates that the interface position $s(t)$ converges to the setpoint $s_r$ monotonically and smoothly without overshooting, i.e., $\dot s(t)>0$ and $s_0<s(t)<s_r$ hold for all $t \geq 0$. Fig. \ref{fig:response} (b) shows that the proposed sampled-data control law maintains constant positive value for every sampling period and is monotonically decreasing to zero. Fig \ref{fig:response} (c) illustrates that the boundary temperature $T(0,t)$ keeps greater than the melting temperature $T_{\rm m}$ with accompanying ``spikes" at every sampling time $t = \tau_{j}$ up to 2 hours. Such spikes are caused by the large drop of the control input $q_{\rm c}(t)$ at sampling time observed from Fig. \ref{fig:response} (b), which affects the boundary temperature directly as given in the boundary condition \eqref{sys2}.  Therefore, the numerical results are consistent with the theoretical results we have established in Lemmas \ref{lem2} and \ref{lem3} for the required properties and in Theorem \ref{theo-1} for the stability of the closed-loop system. 

\section{Conclusion and Future Work}\label{sec:conclusion} 
This paper presented the sampled-data control for the Stefan problem in both one-phase and two-phase cases by application of ZOH to the nominal continuous boundary feedback control law. We proved that under some explicit conditions on the setpoint position and the control gain with respect to the sampling scheduling, the closed-loop system maintains the required conditions for the model validity and is globally exponentially stable. Numerical simulation illustrated the desired performance of the proposed control law. 

While we focused on the full-state feedback design by assuming the availability of the entire temperature profile at each sampling time as a measured value, for the practical implementation it is significant to design an observer-based output feedback control by reconstructing the temperature profile under the availability of only the boundary temperature measured at each sampling time and utilizing the estimated temperature profile as a feedback form, which will be considered as one of our future works. Another interesting direction is ``quantized control" which has a finite or regularly distributed discrete sets of the input value in addition to the sampling time as a \emph{digital} nature \cite{Hayakawa09, Selivanov16}. Owing to the practical implementability of the sampled-data design, a physical experiment of the proposed control law will be demonstrated using some phase change materials.

\bibliographystyle{unsrt}

\end{document}